\documentclass[11pt]{article}

\usepackage{latexsym}
\usepackage{amssymb}
\usepackage{graphicx}
\usepackage{amsmath}
\usepackage{epsfig}
\usepackage{multicol}

\newtheorem{theorem}{Theorem}[section]

\newtheorem{prop}[theorem]{\bf{Proposition}}

\newtheorem{remark}[theorem]{\bf{Remark}}

\newcommand{\conv}{\mathcal{C}}
\newcommand{\V}{\mathcal{V}}
\newcommand{\W}{\mathcal{W}}

\begin{document}

\title{On the Convex Hulls of Self-Affine Fractals}

\author{I. Kirat and I. Kocyigit  }

%\address{Ibrahim Kirat: Department of Mathematics, Istanbul Technical
%University, 34469,Maslak-Istanbul, Turkey }
% \email{ibkst@@yahoo.com}

%%\author{Ibrahim Kirat$^{}$ \\
%%%Ilker Kocyigit$^{b}$ \\
%\\
%$^{a}$
%\footnotesize{\emph{Department of Mathematics, Istanbul Technical University, 34469 Maslak, Istanbul, Turkey}} \\
%\emph{\footnotesize{E-mail: ibkst@yahoo.com}} \\
%$^{b}$\footnotesize{\emph{Department of Mathematics, University of Michigan, Ann Arbor, MI  USA}} \\
%\emph{\footnotesize{E-mail: ilkerk@gmail.com}}
%}

\date { }
%\date { Updated in May 2006 and in January 2010}
\maketitle

\footnote{This unpublished version was written in March 2013.}

\begin{abstract}
Suppose that the set
${\mathcal{T}}= \{T_1, T_2,...,T_q \} $ of real $n\times n$ matrices
has joint spectral radius less than $1$. Then for any digit set $ D=
\{d_1, \cdots, d_q\} \subset {\Bbb R}^n$, there exists a unique
nonempty compact set $F=F({\mathcal{T}},D)$ satisfying $ F = \bigcup
_{j =1}^q T_j(F + d_j)$, which is called a self-affine fractal.
We consider an existing criterion for
the convex hull of $F$ to be a polytope, which is due to  Kirat
and Kocyigit.  %In such a case, we were able
%to find the coordinates and the number of the vertices of the polytope.
%$F$ is called a \emph{finite-type self-affine fractal} if its neighbor
%graph is finite.
In this note, we strengthen our criterion for the case $T_1=T_2=\cdots =T_q $. More
specifically, we give an  upper bound for the number of steps needed for
deciding whether the convex hull of $F$ is a polytope or not.
%stopping-time algorithm to determine whether the convex hull of $F$ is a polytope or not.
This improves our earlier result on the topic.

\end{abstract}

\emph{Keywords} : Self-affine fractal, Polytope convex hulls.

\section{Introduction}\label{INTRO}
As usual, we use ${\Bbb R}$ for the set of real numbers.
Let $M_n({\Bbb R})$ denote the set of $n\times n$ matrices
with real entries. A matrix $T\in M_n({\Bbb R})$ is called  {\it
expanding } (or {\it expansive }) if all its eigenvalues have moduli
$>1$. For a set $\{T_1^{-1},... , T_q^{-1} \}\subset  M_n({\Bbb
R})$ of expanding matrices, the $T_i$ may not be contractive
 with respect to  the same norm (obviously, for each of them there is
a norm, i.e. the Lind norm, depending on the matrix,  with respect to which it is contractive \cite{LW}).

 Let $\|\cdot \|$ denote a
norm on ${\Bbb R}^n$, and $
\| T \|=\sup\{ \| Tx \|: x\in {\Bbb R}^n , \ \ {\| x \|=1} \}
$ be the induced matrix norm of a matrix $T\in M_n({\Bbb R})$. For $\mathcal{T} = \{T_1, \cdots \, T_q \}\subset M_n({\Bbb R})$, we
set
$$
\| \mathcal{T} \|=\max \{ \| T_j \|: \ T_j \in \mathcal{T} \}.
$$
Let $J_k=\{(j_1,j_2,...,j_k): 1\leq j_i \leq q \}$,
and let
${\bf j}= (j_1, \cdots , j_k)$ denote a multi-index or an element in $J_k$
so that $|{\bf j}|=k$ is the length of ${\bf j}$. By $T_{\bf j}$, we mean the product $T_{j_1}\cdots
T_{j_k}$. We let $\mathcal{T}^{k}=\{T_{\bf j}: \ {\bf j}\in J_k  \}.$ Then
the number
$$
\lambda(\mathcal{T})=\lim_k \|\mathcal{T}^{k}\|^{{1}/{k}}=\limsup
\|\mathcal{T}^{k}\|^{{1}/{k}}= \inf_k\|\mathcal{T}^{k}\|^{{1}/{k}}
$$
is called the (uniform) \emph{joint spectral radius} of $\mathcal{T}.$ Then we have the following.
%\begin{defn} The number
%$$
%\lambda(\mathcal{T})=\lim_k \|\mathcal{T}^{k}\|^{{1}/{k}}=\limsup
%\|\mathcal{T}^{k}\|^{{1}/{k}}= \inf_k\|\mathcal{T}^{k}\|^{{1}/{k}}
%$$
%is called the (uniform) joint spectral radius of $\mathcal{T}.$
%\end{defn}

\begin{prop}\label{prop1} {\rm \cite{KK}, \cite{F}} Suppose that
${\mathcal{T}}= \{T_1, T_2,...,T_q \} \subset  M_n({\Bbb R})$
satisfies $\lambda(\mathcal{T}) < 1$. Then for any set $ D= \{d_1,
\cdots, d_q\} \subset {\Bbb R}^n$, called a digit set, there exists a unique nonempty
compact set $F$ satisfying
\begin{equation}\label{eqn1}
F = \bigcup _{j =1}^q T_j(F + d_j).
\end{equation}
%$F = \bigcup _{j =1}^q T_j(F + d_j)$.
\end{prop}

We sometimes write $F(\mathcal{T}, D)$
for $F$ to stress the dependence on $\mathcal{T}$ and $D$. The compact set $F$ in (\ref{eqn1}) is called a {\it self-affine set or a self-affine fractal}, and
can be viewed as the \emph{invariant set} or the attractor of the (affine)
{\it iterated function system} (IFS)
 $\{\phi_j(x)={T_j}(x+d_j)\}_{j=1}^q$. Our generalization of $F$ here includes all classical
self-affine attractors. In the most classical case, $T_j$ are
assumed to be contractions.

 Further, if ${\mathcal{T}}$ is a set of nonsingular matrices with  $\lambda(\mathcal{T}) < 1$, then $\{T_1^{-1},... , T_q^{-1} \}\subset  M_n({\Bbb
R})$ is a set of expanding matrices. In fact,
when $\lambda(\mathcal{T})=1$,
there is no vector norm on ${\Bbb R}^n$ such that the induced matrix norm satisfies
$\| T_j\|<1$ for all $j\in \{1,...,q\}$ (see \cite{KK} for details).

In this paper, we consider $F=F({\mathcal{T}},D)$
in Proposition \ref{prop1} with the assumption that
$ T_1=T_2=\cdots =T_q =T$,
%${\mathcal{T}}\subset  M_n({\Bbb
%R})$
and $ T$ is a nonsingular matrix. The purpose of this note is to report a strengthened
criterion for the convex hull of $F$ to be a polytope. This problem arises from the interest
in the geometry of fractals.

\section{Preliminaries }\label{notation}

In the sequel, we will also assume that $d_1=0$. For ${\bf j}=(j_1,\ldots,j_k)\in J_k$, we let
$\phi_{\bf j}=\phi_{j_1}\circ \cdots \circ \phi_{j_k}$. We set
$$A_k=\{T_{j_1}T_{j_2}\cdots T_{j_k}d_{j_k}+\ldots +
T_{j_1}T_{j_2}d_{j_2}+T_{j_1}d_{j_1} \ \ | \ \ d_{j_1},\ldots,d_{j_k}\in
D\}.$$
%Let ${\mathcal{C}}_n$ denote the space of all non-empty compact
%subsets of ${\Bbb R}^n$. Let $||\cdot ||$ be a norm on ${\Bbb R}^n$.
%The {\it Hausdorff metric} on ${\mathcal{C}}_n$ with respect to
%$||\cdot ||$ is defined by
%
%\begin{equation*}
%d_H(F,F'):=max \{ \sup\limits_{x\in F} \inf\limits_{x'\in F'}
%\left\|x-x'\right\|, \sup\limits_{y'\in F'} \inf\limits_{y\in
%F}\left\|y-y'\right\| \}
%\end{equation*}
%In fact, under the hypothesis of Proposition \ref{prop2}, for any
%initial set $F_0\in {\mathcal{C}}_n$, $\bigcup _{|J|= k}
%\phi_J(F_0)$ converges to $F$ in the Hausdorff metric.

%Let  $|J|=k$ denote the length of a multi-index $J= (j_1, \cdots , j_k)$.
In the Hausdorff metric, we know that
$$A_k=\bigcup_{|{\bf j}|=k} \phi_{\bf j}(0),~~~~~~k=1,2, \ldots, $$
converges to $F$ as $k\rightarrow \infty.$
%In the sequel, we assume
%that $0\in D$.
The {\it convex
hull} of a set $S\subset \mathbb{R}^n$, denoted by $\conv(S)$, is
the intersection of all convex sets in $\mathbb{R}^n$ containing
$S$. It is of geometrical interest to study $\conv(F)$
%\cite{SW1}
and to determine its vertices \cite{SW1}.
Since $0\in D$, it can be easily seen that $A_k\subseteq
A_{k^\prime}$ and $\conv(A_k) \subseteq \conv(A_{k^\prime})$ for
$k\leq k^\prime$. We note that $\conv(F)$
%of a self-affine set $F$
is a compact set.

For brevity, a point $x\in F$ will be
denoted by a sequence $d_{j_1} d_{j_2} \ldots d_{j_k}\ldots.$   A {\it periodic sequence} is a
sequence of the form
$$
%Y=(d_1,d_2,\ldots,d_p,d_1,d_2,\ldots,d_p,d_1,d_2,\ldots)=(
 \overline{d_{j_1}d_{j_2} \ldots d_{j_p}}
 %)
 ,$$ i.e., the block of
digits $d_{j_1},d_{j_2},\ldots,d_{j_p}$ is repeated indefinitely.
 An {\it eventually periodic
(e.p.) sequence} $Y$ is a sequence of the form $$
%Y=(a_1,a_2,\cdots
%a_m, d_1,d_2,\ldots,d_p,d_1,d_2,\ldots,d_p,\ldots)=(
 d_{i_1}d_{i_2}\ldots d_{i_m} \overline{ d_{j_1}d_{j_2} \ldots d_{j_p}}
 %)
 ,$$
 where $ i_1i_2\ldots i_m\in \{1,...,q\}$. For
% bounded $s$ and
all arbitrarily large $r\geq 2$, a point of the form
$$x=i_1i_2\ldots i_m (j_1j_2 \ldots j_p)_rj_{p+1}\ldots j_{p+s}$$ will
be called a  {\it finitely eventually periodic
(f.e.p.) point}. For instance, suppose that $i_1,...,i_m,j_{1},...,j_{p+s}\in \{1,2\}$,
%$s\leq 4$,
then  $(12)_{10}$, $1(12)_{20}112$, $11(12)_{40}1121$, $(12)_{100}21$ are all f.e.p. points.
 Assume that $\W_k\subseteq \V_k$, and $\W_k$ consists of f.e.p. points. Let $\W_{k, \infty}$ be the set of e.p. points of the form
$i_1i_2\ldots i_m \overline{j_1j_2 \ldots j_p}$ associated to the f.e.p.
points of $\W_k\subseteq \V_k$.  Set \\
$$\W_{k, \infty}^C=\{ x=i'_1i_1i_2\ldots i_m
\overline{j_1j_2 \ldots j_p} \ \ | \   \ i_1\ldots i_m
\overline{j_1j_2 \ldots j_p}\in \W_{k, \infty}, \   1\leq i'_1\leq q, \  x\not\in \W_{k,
\infty} \}.$$
By a \emph{string} of a sequence, we mean a special
%of a sequence is finite
sequence consisting of certain consecutive terms of that sequence. For
example, let $Y=d_{j_1}d_{j_2}\ldots d_{j_k}\ldots$ be a sequence of digits, then $d_{j_2} d_{j_3} d_{j_4}$ is a
3-string of $Y$, $d_{j_5} d_{j_6} \ldots d_{j_{k+4}}$ is a k-string of $Y$, and
$d_{j_2} d_{j_3} \ldots$ is an $\infty$-string of $Y$.

In the literature, there are two characterization of self affine-fractals with polytope convex hulls. The first of them is due to Strichartz and Wang, and the other is due to Kirat and Kocyigit . They are as follows:

%\bigskip

\begin{prop}\label{prop2} {\rm \cite{SW1}} Assume that $ T_1=T_2=\cdots =T_q =T$ in (\ref{eqn1}) Let $\{n_j\}$ be the outward unit normal vectors of the
$(n-1)$-dimensional faces of $\conv(D)$. Then $\conv(F)$ is a polytope if and only if
every $n_j$ is an eigenvector of $(T^*)^{-k}$ for some $k$, where $T^*$ is the classical adjoint of $A$.
\end{prop}

%\medskip
%\begin{theorem}\label{single}

\begin{prop}\label{prop3} {\rm \cite{KK}} (i) $\conv(F)$ is a polytope if and only if there exists an index  $k_0$ and
a subset $\W_{k_0}\subseteq \V_{k_0}$ with f.e.p. points such that
%Suppose that there exists a $k_0\in \mathbb{N}$ such that
%\indent
$$\W_{k_0, \infty}^C\cup \V_{k_0} \subseteq \conv(\W_{k_0, \infty}).$$
In such a case, $\conv(F)=\conv(\W_{k_0, \infty}).$

(ii) Assume that $ T_1=T_2=\cdots =T_q =T$ in (\ref{eqn1}). Then $\conv(F)$ is a polytope  if and only if $\#\V_i=\#\V_{i+1}=t$ for some $i$. In such a case,
the points of $\V(F)$ are periodic and for any $k>i$, $\V_{k}$ consists of all $k$-strings of the points of $\conv(F)$.
%\end{theorem}
\end{prop}

%\bigskip

There are two differences between Proposition \ref{prop2} and Proposition \ref{prop3}. First, the former only deals with one matrix. Secondly, unlike the latter, the former doesn't find the coordinates and the number of the vertices of the polytope $\conv(F)$. Although, Proposition \ref{prop3} is better than Proposition \ref{prop2} in that respect, and we can obtain many examples by using it, the second characterization cannot decide if $\conv(F)$ is a polytope when the conditions in Proposition \ref{prop3} are not satisfied for $k_0$ or $i$ up to a large number.

It is our aim in this note to give an upper bound for $i$ satisfying the conditions in Proposition \ref{prop3} in the single-matrix case.
Thus we can say that $\conv(F)$ is not a polytope if the condition in (ii) are not met after a finite number of steps.

%For that, we assume that $F$ is a finite type-fractal and use the concept of neighbor graph of $F$.
%
% Following the terminology of Bandt and Mesing \cite{BM},
%we consider the compositions of the maps $S_1,...,S_q$.
%For ${\bf j}\in J_k$, let $F_{\bf j}=S_{\bf j}(F)=S_{j_1}S_{j_2}\cdots S_{j_k}(F)$.
%If $F_{\bf i}\cap F_{\bf j}\neq \emptyset $ for ${\bf i},{\bf j}\in J_k$, such a pair $F_{\bf i}, F_{\bf j}$ is called
%``a type''. A \emph{neighbor map} has the form $h =S_{\bf i}^{-1}S_{\bf j}$ where
%$F_{\bf i}\cap F_{\bf j}\neq \emptyset $.
%Let these maps form the vertex set of a graph $G$ with root vertex id, the identity map. The edges of the graph lead from
%each vertex $h$ to another vertex $g = S_i^{-1}h S_j$ with $i, j\in \{1,2,...,q \} $ and are labelled with the
%corresponding pair of symbols $i, j$.
%
%The graph G will be called the \emph{neighbor graph} of the fractal F. We call F a \emph{finite-type self-affine fractal} if the graph $G$ is finite (cf. \cite{BM}). In this case, from each vertex there starts a path which ends in a directed cycle. From now on, we only consider finite-type self-affine fractals, which includes all self-affine tiles. The length of a cyclic path is the number of vertices of that path. We restrict our attention to only those cyclic paths which give different periodic sequences.  Let $l_c$ be the sum of the lengths of all cyclic paths in the neighbor graph of $F$. Then we have the following result.

\section{Statement of the Result}\label{strengthened criterion}

${\Bbb C}$ stands for the set of complex numbers. In the following theorem, $z_n\in {\Bbb C}$ denotes an $n$-th root of unity.

\begin{theorem}\label{single criterion}
Assume that $ T_1=T_2=\cdots =T_q =T$ in (\ref{eqn1}). Let
$$U=\{c_{n_1}z_{2n_1}, c_{n_2}z_{2n_2}, \cdots , c_{n_m}z_{2n_m}\}$$ be
the set of all roots of $T^{-1}$ of the form $cz_n$, where $c> 0$. 
%Assume that $U\neq \emptyset$.
Set $k=2n_1n_2\cdots n_m$.
Then if $\#\V_i\neq \#\V_{i+1}$ for all $i\leq k$, then $\conv(F)$ is not a polytope. Therefore, $\conv(F)$ is a polytope
if and only if $\#\V_i= \#\V_{i+1}$ for some $i\leq k$.
\end{theorem}

\begin{remark}\label{}.\rm{ Here we don't assume that $ T^{-1}$ is a similitude, it is a general expanding matrix. }
%It follows from Proposition \ref{prop2} that
%$\conv(F)$ is not a polytope when $U=\emptyset$. Note that $U=\emptyset$ is possible only if $F\subset {\Bbb R}^{2n}$.}
\end{remark}

%
%\medskip
%
%
%\noindent \textbf{Acknowledgments}
%
%
%{\footnotesize Part of this work was done while the first author was visiting the
%Department of Mathematics, the Chinese University of Hong Kong, and Institute for Mathematics, Arndt University, Greifswald, Germany.
%He would like to thank Professor Ka-Sing Lau for his kind invitation, constant support and helpful discussion, and Professor Christoph Bandt for his hospitality and valuable discussion. Thanks are extended to De-Jun Feng, Antti K${\rm \ddot{a}}$enm${\rm \ddot{a}}$ki, Boris Solomyak and Pablo Shmerkin for valuable conversations.}


{\footnotesize \begin{thebibliography}{9999}



\bibitem {F}   K.J. Falconer, Fractal geometry: Mathematical
Foundations and Applications, John Wiley $\&$ Sons, Chichester, 1990.

\bibitem {KK} I. Kirat and I. Kocyigit, Remarks on  self-affine  fractals  with polytope convex hulls, Fractals., 18 (2010) no.4,  pp. 483-498.

\bibitem {LW} \ J. C. Lagarias and  Y. Wang,  Integral self-affine
tiles in ${\Bbb R}^n$, {\it Adv. Math.}, {\textbf{121}} (1996) 21-49.

\bibitem{SW1} R. Strichartz, and  Y. Wang,
Geometry of self-affine tiles I, {\it Indiana Univ. Math. \ J } \textbf{48}
(1999) 1-23.


\vspace{1cm}
\noindent Ibrahim Kirat           \hfill  Ilker Kocyigit \\
Department of Mathematics         \hfill  Department of Mathematics\\
Istanbul Technical University     \hfill  University of Michigan \\
Maslak 34469, Istanbul            \hfill  Ann Arbor, MI 48109-1043     \\
Turkey                            \hfill  USA \\

\noindent E-mail:  ibkst@yahoo.com  \hfill E-mail: ilkerk@gmail.com


\end{thebibliography} }
\end{document}